\RequirePackage{fix-cm}
\documentclass[smallextended]{article}       
\usepackage{graphicx}
\usepackage{amsfonts}
\usepackage{amssymb}
\usepackage[all]{xy}
\newcommand{\id}{\mathop{\rm id}\nolimits}
\newcommand{\End}{\mathop{\rm End}\nolimits}
\newcommand{\phe}{\varphi}
\newcommand{\C}{\mathbb{C}}
\newcommand{\Q}{\mathbb{Q}}

\newcommand{\N}{\mathbb{N}}
\newcommand{\Z}{\mathbb{Z}}
\newcommand{\PP}{\mathbb{P}}
\newtheorem{definition}{Definition}{\bfseries}{\rmfamily}
\newtheorem{remark}{Remark}{\itshape}{\rmfamily}
 \newtheorem{theorem}{Theorem}{\bfseries}{\itshape}
  \newtheorem{corollary}{Corollary}{\bfseries}{\itshape}
   \def\ackname{Acknowledgements}
  \def\acknowledgement{\par\addvspace{17pt}\small\rmfamily
\trivlist\if!\ackname!\item[]\else
\item[\hskip\labelsep
{\bfseries\ackname}]\fi}

\newenvironment{acknowledgements}{\begin{acknowledgement}}
{\end{acknowledgement}}
\newcommand{\eqref}[1]{(\ref{#1})}
\newcounter{subesercizio}
\newenvironment{thlist}%
{\begin{list}{\rm(\roman{subesercizio})}{\usecounter{subesercizio}%
\setlength{\labelsep}{3pt}\setlength{\leftmargin}{16pt}\setlength{\topsep}{3pt}%
\addtolength{\itemsep}{1pt}%
\setlength{\labelwidth}{\leftmargin}%
\setlength{\listparindent}{\itemindent}}}{\end{list}}

\begin{document}

\title{Open problems in local discrete holomorphic dynamics
}


\author{Marco Abate}

\maketitle



\maketitle

\begin{abstract}
This paper contains a selection, dictated by personal taste and by no means complete, of open
problems in local discrete holomorphic dynamics.
\smallskip


\end{abstract}

\section{Introduction}
\label{intro}
The aim of this paper is to collect and put in context a few open problems in the 
area of local holomorphic discrete dynamics; so let us begin by defining what 
is a discrete local holomorphic dynamical system.

\begin{definition}
Let $M$ be a complex manifold, and $p\in M$. A \emph{discrete holomorphic local dynamical system} at~$p$ is a holomorphic map
$f\colon U\to M$ such that~$f(p)=p$, where $U\subseteq M$ is an open
neighbourhood of~$p$; we shall also always assume that $f\not\equiv\id_U$. We shall
denote by~$\End(M,p)$ the set of holomorphic local dynamical systems at~$p$. 
\end{definition}

We shall be mainly concerned with the behavior of~$f$
nearby~$p$, and thus $\End(M,p)$ actually is the set of germs of holomorphic
self-maps of~$M$ at~$p$; for this reason we shall often use the word ``germ" as an abbreviation of
``discrete holomorphic local dynamical system", and we shall
allow us to restrict the domain of $f\in\End(M,p)$ to any suitable
neighbourhood of~$p$ whenever useful.

%
In the next sections we shall discuss some specific open problems, selected according to our taste and with no pretense of completeness; but in this introduction
we shall instead present the most basic questions one can ask about a local discrete
dynamical system. Claiming to have completely understood a given local dynamical systems 
amounts to having complete answers to at least the first two questions below; and
claiming to have understood a given \emph{class} of discrete holomorphic local dynamical systems
most of the times amounts to having answers to the remaining three questions.

First of all, a local dynamical system is not a dynamical system in the standard (global) sense
of the word, because points can escape from the domain of definition. However, it is easy
to associate a \emph{bona-fide} dynamical system with any local dynamical system.
The phase space of this new dynamical system is the stable set:

\begin{definition}
Let $f\in\End(M,p)$ be a (discrete holomorphic) local dynamical system defined
on an open set $U\subseteq M$. Then the \emph{stable set}~$K_f$ of~$f$ is
$$
K_f=\bigcap_{k=0}^\infty f^{-k}(U)\;.
$$
In other words, the stable set of~$f$ is
the set of all points~$z\in U$ such that the \emph{orbit} $\{f^k(z)\mid
k\in\mathbb{N}\}$ is well-defined, where $f^k$ denotes the $k$-th iterate of~$f$. If $z\in U\setminus K_f$, we shall say that $z$ (or
its orbit) {\sl escapes} from~$U$. 
\end{definition}

Clearly, $p\in K_f$, and so the stable set is never empty (but it can happen
that $K_f=\{p\}$). It depends a priori on~$U$, but again in most cases we shall be
interested only in the behavior nearby~$p$. 
Thus the first natural question in discrete holomorphic local dynamics is:
\begin{itemize}
\item[(Q1)] \emph{What is the topological structure of (the germ at~$p$ of)~$K_f$?}
\end{itemize}

\noindent For instance, does $K_f$ have non-empty interior? Is it locally connected at~$p$?
Is $K_f\setminus\{p\}$ connected? What is the topological (homological, cohomological) structure
of~$U\setminus K_f$? And so on.

\begin{remark}
Both the definition of stable set and Question~1 are topological in character;
we might also state them for local dynamical systems which are continuous only. However, the {\it answers} might (and usually will) strongly depend on the
holomorphicity of the dynamical system.
\end{remark}


Clearly, the stable set~$K_f$ is completely $f$-invariant, and thus
the pair $(K_f,f)$ is the promised discrete global dynamical system, canonically associated with the given discrete local dynamical system. In particular, the second natural question in discrete holomorphic local dynamics is
\begin{itemize}
\item[(Q2)] \emph{What is the dynamical structure of~$(K_f,f)$?}
\end{itemize}
\noindent For instance, what is the asymptotic behavior of the orbits? Do
they converge to~$p$, or have they a chaotic behavior? Is there a dense
orbit? Do there exist proper $f$-invariant subsets, that is sets
$L\subset K_f$ such that~$f(L)\subseteq L$? If they do exist, what is the
dynamics on them?

To answer all these questions, one of the most efficient ways is to replace $f$ by a
dynamically equivalent but simpler (e.g., linear) system~$g$. In our context,
``dynamically equivalent" means ``locally conjugated"; and we have different 
kinds of conjugacy to consider. 

\begin{definition}
Let $f_1\colon U_1\to M_1$ and $f_2\colon U_2\to M_2$ be two holomorphic local
dynamical systems at~$p_1\in M_1$ and~$p_2\in M_2$ respectively. We shall say
that~$f_1$ and~$f_2$ are \emph{holomorphically} (respectively, \emph{smoothly,} $C^k$ with
$k\in\mathbb{N}^*$, or \emph{topologically}) \emph{locally conjugated} if there are open neighbourhoods
$W_1\subseteq U_1$ of~$p_1$, $W_2\subseteq U_2$ of~$p_2$, and a
biholomorphism (respectively, a $C^\infty$ diffeomorphism, a $C^k$ diffeomorphism, or a homeomorphism) $\phe\colon W_1\to W_2$ with
$\phe(p_1)=p_2$ such that 
$$
f_1=\phe^{-1}\circ f_2\circ\phe
$$
on $\phe^{-1}\bigl(
W_2\cap f_2^{-1}(W_2)\bigr)=W_1\cap f_1^{-1}(W_1)$.  
\end{definition}

If $f_1\colon U_1\to M_1$ and $f_2\colon U_2\to M_2$ are locally conjugated we clearly have
$$
K_{f_2|_{W_2}}=\phe(K_{f_1|_{W_1}})\;;
$$
so the local dynamics of~$f_1$
about~$p_1$ is to all purposes equivalent (up to the order of smoothness of the local
conjugation) to the local dynamics of~$f_2$
about~$p_2$.

In particular,  using local coordinates centered at~$p\in M$ it is easy to show that any
holomorphic local dynamical system at~$p$ is holomorphically locally conjugated
to a holomorphic local dynamical system at~$O\in\mathbb{C}^n$, where $n=\dim M$;
so from now on we shall mostly work only with $\End(\mathbb{C}^n,O)$.

Whenever we have an equivalence relation in a class of objects, 
classification problems come out. So the third natural question in local
holomorphic dynamics is
\begin{itemize}
\item[(Q3)] \emph{Find a (possibly small) class $\mathcal{F}$ of holomorphic
    local dynamical systems at~$O\in\mathbb{C}^n$ such that every
    holomorphic local dynamical system~$f$ at a point in an
    $n$-dimensional complex manifold is holomorphically (respectively,
    smoothly, $C^k$, or topologically) locally conjugated to a (possibly) unique element
    of~$\mathcal{F}$, called \emph{holomorphic} (respectively, \emph{smooth,} $C^k$ or
      \emph{topological}) \emph{normal form} of~$f$.} 
\end{itemize}
Unfortunately, the holomorphic classification is often too
complicated to be practical; the family~$\mathcal{F}$ of normal forms might
be uncountable. A possible replacement is looking for invariants
instead of normal forms: \smallskip
\begin{itemize}
\item[(Q4)] \emph{Find a way to associate a (possibly small) class of
    (possibly computable) objects, called \emph{invariants,} to any
    holomorphic local dynamical system~$f$ at~$O\in\mathbb{C}^n$ so
    that two holomorphic local dynamical systems at~$O$ can be
    holomorphically (respectively, smoothly, $C^k$, or topologically) locally conjugated only if they have the same invariants.
    The class of invariants is \emph{complete} if two
    holomorphic local dynamical systems at~$O$ are holomorphically (respectively, smoothly, $C^k$, or topologically) locally
    conjugated if and only if they have the same invariants.}
\end{itemize}

\noindent Contrarily to the previous ones, our final general question makes sense only
for \emph{holomorphic} local dynamical systems. Indeed,
a discrete holomorphic local dynamical system at~$O\in\mathbb{C}^n$ is 
given by an element of~$\mathbb{C}_0\{z_1,\ldots,z_n\}^n$, the space
of $n$-tuples of converging power series in~$z_1,\ldots,z_n$ without
constant terms, which is a
subspace of the space~$\mathbb{C}_0[\![z_1,\ldots,z_n]\!]^n$ of $n$-tuples
of formal power series without constant terms. It is well known that 
$\mathbb{C}_0[\![z_1,\ldots,z_n]\!]^n$ is closed under composition of power series, 
and that an element $\Phi\in
\mathbb{C}_0[\![z_1,\ldots,z_n]\!]^n$ has an inverse (with respect to
composition) still belonging to~$\mathbb{C}_0[\![z_1,\ldots,z_n]\!]^n$ if
and only if its linear part is a linear automorphism
of~$\mathbb{C}^n$.

\begin{definition}
 We 
say that 
$f_1$,~$f_2\in\mathbb{C}_0[\negthinspace[z_1,\ldots,z_n]\negthinspace ]^n$ are \emph{formally conjugated} if there
is an invertible $\Phi\in\mathbb{C}_0[\![z_1,\ldots,z_n]\!]^n$ such that
$f_1=\Phi^{-1}\circ f_2\circ\Phi$ in~$\mathbb{C}_0[\![z_1,\ldots,z_n]\!]^n$. 
\end{definition}

Clearly, two holomorphically locally conjugated holomorphic
local dynamical systems are both formally and topologically locally
conjugated too. On the other hand, there are examples of
holomorphic local dynamical systems that are topologically locally
conjugated without being neither formally nor holomorphically locally
conjugated, and examples of holomorphic local dynamical systems that
are formally conjugated without being neither holomorphically nor
topologically locally conjugated. So the last natural general question in
local holomorphic dynamics is
\begin{itemize}
\item[(Q5)] \emph{Find normal forms and invariants with respect to the relation
of formal conjugacy for holomorphic local dynamical systems at~$O\in\mathbb{C}^n$.}
\end{itemize}

In the rest of this paper we shall describe a few specific open problems, 
both in one and in several complex variables; we
refer to \cite{A5,A6} and \cite{Br2} for surveys on the theory of discrete holomorphic local dynamical systems, and for more details on the known answers to the previous questions and on background of the open problems we selected. 

\section{One complex variable}
\label{sec:1}
\subsection{Linearization in specific families}
\label{sec:1.1}

A discrete holomorphic local dynamical system in one complex variable is given by a germ of holomorphic
function fixing the origin of the form
\[
f(z)=\lambda z+a_2 z^2+\cdots\in\C\{z\}\;,
\]
where $\lambda=f'(0)\in\C^n$ is the \emph{multiplier} of $f$. It is well-known (K\oe nigs' theorem) that
if $|\lambda|\ne 0$,~$1$ then $f$ is \emph{holomorphically linerizable,} that is locally holomorphically conjugated to the linear map $w\mapsto \lambda w$ --- and thus in this case questions (Q1)--(Q5) presented in the
introduction are easily solved. If $a_1=0$, and then we can write $f(z)=a_k z^k+o(z^k)$ with
$k\ge 2$ and $a_k\ne 0$, it is also well-known (B\"ottcher's theorem) that $f$ is locally holomorphically conjugated to
$w\mapsto w^k$, and thus in this case too the local dynamics is completely clear.

If $|\lambda|=1$ and $\lambda=e^{2\pi i p/q}$ is a $q$-th root of unity, then it is easy to prove that 
$f$ is (topologically, holomorphically or formally) linearizable if and only if $f^q=\id$; therefore the \emph{linearization problem}, that is
deciding when $f$ is linearizable, is meaningful only when $\lambda=e^{2\pi i\theta}$ with $\theta\notin\Q$. Note that it is not difficult to show that every germ with multiplier of this form is formally
linearizable, and that it is topologically linearizable if and only if it is holomorphically linearizable;
so the main question here is deciding whether a given germ is holomorphically linearizable
or not --- and, in particular, whether this question can be solved just by examining the multiplier.

The main result in this area is due to Brjuno \cite{Brj1,Brj2,Brj3} and Yoccoz \cite{Y1,Y2}. To state it, let us
introduce a bit of terminology and some notations. 

\begin{definition}
\label{def:1.1}
We shall say that $f\in\End(\C,O)$ is \emph{elliptic} if its multiplier has modulus one but
is not a root of unity; and that
the origin is a \emph{Siegel point} (respectively,
a \emph{Cremer point}) if $f$ is (respectively, is not) holomorphically linearizable.
\end{definition}

\begin{definition}
\label{def:1.2}
For $\lambda\in S^1$ and $m\ge 1$ put
\[
\Omega_\lambda(m)=\min_{1\le k\le m}|\lambda^k-\lambda|\;.
\]
Clearly, $\lambda$ is a root of unity (or $\lambda=0$) if and
only if $\Omega_\lambda(m)=0$ for all $m$ greater than or equal to some~$m_0\ge 1$;
furthermore, if $|\lambda|\ne 0$,~1 then $\Omega_\lambda(m)$ is bounded away from zero, whereas
if $|\lambda|=1$ then
\[
\lim_{m\to+\infty}\Omega_\lambda(m)=0\;.
\] 
We shall say that $\lambda\in S^1$, not a root of unity, satisfies the \emph{Brjuno condition} (or that $\lambda$ is a \emph{Brjuno number}) if
\begin{equation}
\sum_{k=0}^{+\infty}{1\over2^k}\log{1\over\Omega_\lambda(2^{k+1})}<+\infty\;.  
  \label{eqqsei}
\end{equation}
\end{definition}

\begin{remark}
\label{rem:1.1}
There are several equivalent reformulations of the Brjuno condition, as the convergence of other series involving $\Omega_\lambda$, or as the convergence of a series involving the continuous fraction expansion of~$\theta$; see, e.g., \cite{Rus},
\cite{R5} and \cite{Y2}.
\end{remark}

We can then state the famous \emph{Brjuno-Yoccoz theorem:}

\begin{theorem}[Brjuno, 1965 \cite{Brj1,Brj2,Brj3}, Yoccoz, 1988 \cite{Y1,Y2}]
\label{th:1.BY}  
 Let
$\lambda\in S^1$, not a root of unity. Then the following statements are equivalent:
\begin{thlist}
\item the origin is a Siegel point for the quadratic polynomial $f_\lambda(z)=\lambda
z+z^2$; 
\item the origin is a Siegel point for all 
$f\in\End(\mathbb{C},0)$ with multiplier~$\lambda$;
\item the number $\lambda$ satisfies the Brjuno condition.
\end{thlist}
\end{theorem}

This theorem has two aspects. From one side, it gives a necessary and sufficient condition
on $\lambda$ ensuring that \emph{local} dynamical systems with multiplier $\lambda$ are
holomorphically linearizable. However, there clearly exist holomorphically linearizable germs
whose multiplier does not satisfy the Brjuno condition: it suffices to take a germ of the form
$f(z)=\phe^{-1}\bigl(\lambda\phe(z)\bigr)$, where $\phe$ is a local biholomorphism and
$\lambda$ is not a Brjuno number. 

One the other hand, Theorem~\ref{th:1.BY} says that in the family of quadratic polynomials
a given germ is holomorphically linearizable \emph{if and only if} its multiplier is a Brjuno number.
This observation immediately leads to the first open problem of this survey:
\begin{itemize}
\item[(OP1)] \emph{Find families $\{f_{\lambda, a}(z)=\lambda z+z^2 g_a(z)\}_{a\in M}\subset\End(\C,O)$ of elliptic discrete holomorphic
dynamical systems such that $f_{\lambda,a}$ is holomorphically linearizable if and only
if $\lambda$ is a Brjuno number.} For instance, \emph{is it true that an elliptic polynomial is holomorphically
linearizable if and only if its multiplier is a Brjuno number?} Or, even more specifically, 
\emph{is it true that an elliptic cubic polynomial is holomorphically linearizable if and only if its multiplier
is a Brjuno number?}
\end{itemize}

In this context, it might be useful to remember the following dichotomy due to Ilyashenko and Perez-Marco:

\begin{theorem}[Ilyashenko, 1979 \cite{I1}, Perez-Marco, 2001 \cite{P9}]
\label{th:1.IPM}
Given $g\in\End(\C,O)$ and $\lambda\in \C^*$ not a root of unity, put
 \[
 f_{\lambda,a}=\lambda z+azg(z)
 \]
 for all $a\in\C$. Then:
 \begin{itemize}
 \item[\textup{(i)}] either the origin is a Siegel point of $f_{\lambda,a}$ for all $a\in\C$, or
 \item[\textup{(ii)}] the origin is a Cremer point of $f_{\lambda,a}$ for all $a\in\C\setminus K$, where
 $K\subset\subset\C$ is a bounded exceptional set of capacity (and hence Lebesgue measure) zero.
 \end{itemize}
 \end{theorem}

\subsection{Regularity of the Brjuno function}
\label{sec:1.2}

There is a different way of expressing the Brjuno condition, leading to another interesting open problem.

Given $\theta\in[0,1)$ set
\begin{eqnarray*}
r(\theta)=\inf\{r(f)&\mid& f\in\End(\mathbb{C},0)\hbox{ is
defined and injective in $\Delta$}\\
&&\hbox{and has multiplier $e^{2\pi i\theta}$}\}\;,
\end{eqnarray*}
where $\Delta\subset\C$ is the unit disk and $r(f)\ge 0$ is the radius of convergence of the unique formal linearization of~$f$
with multiplier~1.

On the other hand, given an irrational number $\theta\in[0,1)$ let $\{p_k/q_k\}$ be
the sequence of rational numbers converging to~$\theta$ given by the expansion in continued
fractions, and put
\begin{eqnarray*}
\alpha_n&=&-{q_n\theta-p_n\over q_{n-1}\theta-p_{n-1}},\qquad \alpha_0=\theta,\\
\beta_n&=&(-1)^n(q_n\theta-p_n),\qquad\beta_{-1}=1.  
\end{eqnarray*}

\begin{definition}
\label{def:1.3}
  The {\sl Brjuno function} $B\colon[0,1)\setminus\mathbb{Q}\to(0,+\infty]$ is defined by
$$
B(\theta)=\sum_{n=0}^\infty \beta_{n-1}\log{1\over\alpha_n}.
$$
\end{definition}

Then Yoccoz has proved the following quantitative relationship between the 
infimum $r(\theta)$ of the radii of convergence and the Brjuno function:

\begin{theorem}[Yoccoz, 1988 \cite{Y2}]
\label{th:1.Yoccozest}
\begin{thlist}
\item $B(\theta)<+\infty$ if and only if
$\lambda=e^{2\pi i\theta}$ is a Brjuno number;
\item there exists a universal constant $C>0$ such that
$$
|\log r(\theta)+B(\theta)|\le C
$$
\indent for all $\theta\in[0,1)\setminus\mathbb{Q}$ such that $B(\theta)<+\infty$;
\item if $B(\theta)=+\infty$ then there exists a
  non-linearizable $f\in\End(\mathbb{C},0)$ with multiplier~$e^{2\pi
    i\theta}$.
\end{thlist}
\end{theorem}

The Brjuno function is clearly quite an irregular function, diverging at a $G_\delta$-dense set
of irrational numbers . A surprising fact is that, on the contrary, 
numerical experiments (see, e.g., \cite{Car} and references therein) as well as theoretical
results (see \cite{MMY})
suggest that $\log r+B$ is more regular, even better than continuous. More specifically, 
Marmi, Mattei and Yoccoz proposed the following problem:

\begin{itemize}
\item[(OP2)] \emph{Does $\log r+B$ admit a $1/2$-H\"older
continuous extension to~$[0,1)$?}
\end{itemize}

\subsection{Classification of Cremer points}
\label{sec:1.3}

The local dynamics about a Siegel point is completely clear. The local dynamics about a Cremer point, on the other hand,
is extremely complicated. The best results up to now are due to Perez-Marco and Biswas:

\begin{theorem}[P\'erez-Marco, 1995 \cite{P6,P7}]
\label{th:1.PerezMarco}
Assume that $0$
is a Cremer point for an elliptic discrete holomorphic local dynamical system
$f\in\End(\mathbb{C},0)$. Then:
\begin{thlist}
\item The stable set $K_f$ is compact, connected, full
  (i.e., $\mathbb{C}\setminus K_f$ is connected), it is not reduced
  to~$\{0\}$, and it is not locally connected at any point distinct
  from the origin.
\item Any point of~$K_f\setminus\{0\}$ is recurrent (that is, 
a limit point of its orbit).
\item There is an orbit in~$K_f$ which accumulates at the
  origin, but no non-trivial orbit converges to the origin.
\end{thlist}
\end{theorem}

\begin{theorem}[Biswas, 2010 \cite{B2}]
\label{th:1.Biswas}
The multiplier and the conformal class of the stable set $K_f$ are a complete
set of holomorphic invariants for Cremer points. In other words, two
elliptic non-linearizable holomorphic local dynamical systems 
are holomorphically locally conjugated if and only if they have the
same multiplier and there is a biholomorphism (not necessarily conjugating the dynamics) of a neighbourhood
of~$K_f$ with a neighbourhood of~$K_g$.
\end{theorem}

Surprisingly enough, the topological classification of Cremer points is still open. Clearly, the
homeomorphism class of the (germ at the origin of the) stable set is a topological invariant;
moreover, a non-trivial theorem due to Nashul (see \cite{P7} for another proof) shows that
the multiplier is another topological invariant:

\begin{theorem}[Naishul, 1983 \cite{N}]
\label{th:1.Naishul}
Let $f$, $g\in\End(\mathbb{C},O)$ be two elliptic discrete holomorphic local dynamical systems. If $f$ and $g$ are topologically locally conjugated
then $f'(0)=g'(0)$.
\end{theorem}

Thus a natural open question in this context is:
\begin{itemize}
\item[(OP3)]\emph{Are the multiplier and the homeomorphism class of the stable set a complete set
of topological invariants for discrete holomorphic local dynamical systems in $\End(\C,O)$ having a
Cremer point at the origin?}
\end{itemize}

\subsection{Effective classification of parabolic germs}
\label{sec:1.4}

\begin{definition}
\label{def:1.4}
A local dynamical system $f\in\End(\C,0)$ is \emph{parabolic} if its multiplier is a root of unity;
it is \emph{tangent to the identity} if its multiplier is~1.
\end{definition}

Clearly, if $f\in\End(\C,0)$ is parabolic then a suitable iterate $f^q$ is tangent to the identity; so
most dynamical questions for parabolic systems can be reduced to the study of tangent to the identity germs.

A qualitative description of the dynamics of tangent to the identity germs is given by the
famous \emph{Leau-Fatou flower theorem.} To state it we need to recall a couple of
definitions:

\begin{definition}
\label{def:1.5}
Let $f\in\End(\mathbb{C},0)\setminus\{\id\}$ be tangent to the identity, and thus of the form
\[
f(z)=z+a_{r+1}z^{r+1}+O(z^{r+2}
\] 
with $a_{r+1}\ne 0$, where $r+1\ge2$ is the
\emph{multiplicity} of~$f$. A unit vector~$v\in S^1$
is an \emph{attracting} (respectively, \emph{repelling}) \emph{direction}
for~$f$ at the origin if~$a_{r+1}v^r$ is real and negative (respectively,
positive). 
\end{definition}

Clearly, there are $r$ equally spaced attracting directions,
separated by $r$ equally spaced repelling directions. Furthermore, a repelling
(attracting) direction for~$f$ is attracting (repelling) for~$f^{-1}$, which is
defined in a neighbourhood of the origin.

It turns out that to every attracting direction is associated a connected
component of~$K_f\setminus\{0\}$. 
\begin{definition}
\label{def:1.6}
Let~$v\in S^1$ be an attracting direction for
an~$f\in\End(\mathbb{C},0)$ tangent to the identity. The \emph{basin} centered at~$v$ is the set of
points~$z\in K_f\setminus\{0\}$ such that $f^k(z)\to 0$ and $f^k(z)/|f^k(z)|\to
v$ (notice that, up to shrinking the domain of~$f$, we can assume that $f(z)\ne
0$ for all $z\in K_f\setminus\{0\}$). If $z$ belongs to the basin centered
at~$v$, we shall say that the orbit of~$z$ \emph{tends to~$0$ tangent to~$v$.}

 An \emph{attracting petal} centered at an attracting direction~$v$ of $f$ is an open simply
connected $f$-invariant set $P\subseteq K_f\setminus\{0\}$ such that a point
$z\in K_f\setminus\{0\}$ belongs to the basin centered at~$v$ if and only if its
orbit intersects~$P$. In other words, the orbit of a point tends to~$0$ tangent
to~$v$ if and only if it is eventually contained in~$P$. A \emph{repelling petal}
(centered at a repelling direction) is an attracting petal for the inverse
of~$f$.
\end{definition}

Then:

\begin{theorem}[Leau, 1897 \cite{L}; Fatou, 1920 {[32--34]}] 
 \label{th:1.flower} 
Let
$f\in\End(\mathbb{C},0)$ be a discrete holomorphic local dynamical system tangent to the identity
with multiplicity~$r+1\ge 2$ at the fixed point. Let $v_1^+,\ldots,v_r^+\in
S^1$ be the
$r$ attracting directions of~$f$ at the origin, and $v_1^-,\ldots,v_r^-\in
S^1$ the $r$ repelling directions. Then 
\begin{thlist}
\item for each attracting (repelling) direction~$v_j^\pm$
there exists an attracting (repelling) petal~$P_j^\pm$, so that the
union of these $2r$ petals is a pointed neighbourhood
of the origin. Furthermore, the
$2r$ petals are arranged ciclically so that two petals intersect if and only if
the angle between their central directions is~$\pi/r$.
\item $K_f\setminus\{0\}$ is the (disjoint) union of the basins
centered at the
$r$ attracting directions.
\item If $B$ is a basin centered at one of the attracting directions then
there is a function $\phe\colon B\to\mathbb{C}$ such that \indent $\phe\circ
f(z)=\phe(z)+1$ for all $z\in B$. Furthermore, if $P$ is the corresponding petal
constructed in part {\rm (i),} then $\phe|_P$ is a biholomorphism
with an open subset of the complex plane containing a right
half-plane --- and so $f|_P$ is holomorphically
conjugated to the translation $z\mapsto z+1$.
\end{thlist}
\end{theorem}

Starting from this theorem, Camacho \cite{C} and, independently, Shcherbakov \cite{S} have completed the topological classification of germs tangent to the identity, showing the the multiplicity is a complete set of topological invariants:

\begin{theorem}[Camacho, 1978 \cite{C}; Shcherbakov, 1982 \cite{S}]
\label{th:1.Camacho}
Any discrete hol\-om\-or\-phic local dynamical system tangent to the identity
with multiplicity~$\nu\ge 2$ is topologically locally
conjugated to $g(z)=z+z^\nu$.
\end{theorem}

Furthermore, the formal classification is obtained with not too difficult a computation, and a complete set of invariants is given by the multiplicity
and another complex number, the \emph{index,} explicitly computable. 

On the other hand, the holomorphic classification turned out to be incredibly more difficult.
\'Ecalle \cite{E1,E2,E3,E4} and, independently, Voronin \cite{V} gave a complete set of invariants,
consisting in the multiplicity, the index, and a functional invariant, an equivalence class of functions with specific properties constructed starting from the biholomorphisms introduced in
Theorem~\ref{th:1.flower}.(iii). This set of invariants is not only complete but also full, in the sense
that every possible value of the invariants is realized by a germ tangent to the identity; however, to explicitly
compute \'Ecalle-Voronin functional invariant is an almost impossible task. In particular,
the following problem is still open:
\begin{itemize}
\item[(OP4)] \emph{Give an effective procedure for deciding whether two germs tangent to the
identity are holomorphically locally conjugated.}
\end{itemize}

Clearly, a similar question can be asked for parabolic germs; in that case another (topological, formal and holomorphic) invariant is the multiplier, and one has to replace the multiplicity 
by a suitably defined parabolic multiplicity, and the index by a suitably defined parabolic index
(or, better yet, by \'Ecalle's iterative residue).
Then the multiplier and the parabolic multiplicity are a complete set of topological invariants,
the multiplier, the parabolic multiplicity and the iterative residue are a complete set of
formal invariants, and adding suitably adjusted \'Ecalle-Voronin functional invariants one obtains a
complete set of holomorphic invariants. Again, most of the times the computation of the
functional invariants is hopeless, and thus we have the following generalization of
the previous question:
\begin{itemize}
\item[(OP$4'$)] \emph{Give an effective procedure for deciding whether two parabolic germs are holomorphically locally conjugated.}
\end{itemize}

See also \cite{I1}, \cite{M1,M2} and \cite{MR} for alternative presentations of \'Ecalle-Voronin
invariants.

\section{Several complex variables}
\label{sec:2}

We have seen how in one complex variable the multiplier (that is, the derivative at the fixed point)
plays a fundamental role. In several complex variables instead of the multiplier we may 
consider the (eigenvalues of the) differential at the fixed point, and give a first classification of
discrete holomorphic local dynamical systems based on them. Clearly there are many cases
to consider, and correspondingly many ways to precise the five basic questions we stated in the introduction. Here we shall limit ourselves to a selection of some important open problems
for three main classes of systems: \emph{non-invertible,} \emph{tangent to the identity,} 
and \emph{linearizable} systems.

A piece of terminology we shall systematically use is the following:

\begin{definition}
\label{def:2.he}
The \emph{homogeneous expansion} of a $f\in\End(\C^n,O)$ is the expansion
\[
f(z)=\sum_{j\ge c(f)} P_j(z)
\]
where $c(f)\ge 1$ is the \emph{order} of~$f$, and $P_j$ is an $n$-tuple of homogeneous
polynomials of degree~$j$ (and we are of course assuming that $P_{c(f)}\not\equiv O$). 
Furthermore, we shall say that $f$ is \emph{dominant} if $\det\hbox{Jac}(f)\not\equiv O$. 
\end{definition}

\subsection{Non-invertible systems: asymptotic attraction rate}
\label{sec:2.1}

A discrete holomorphic local dynamical system $f\in\End(\C^n,O)$ is invertible if and only if
the differential $df_O$ is; therefore $f$ is non-invertible if and only if 0 is an eigenvalue of $df_O$. 
In particular, we shall also say that $f$ is \emph{superattracting} if $df_O\equiv O$, i.e., 
if $c(f)\ge 2$.

There are a couple of interesting open questions (suggested by Mattias Jonsson) about
the sequence $\{c(f^k)\}_{k\in\N}$. In dimension~1 it is easy to see that $c(f^k)=c(f)^k$
for all $k\ge 1$, and thus $c(f^k)^{1/k}=c(f)$ for all $k\ge 1$. 
On the other hand, in dimension 2 or more we only have
\begin{equation}
c(f^{h+k})\ge c(f^h)c(f^k)
\label{eq:2.1}
\end{equation}
for all $h$,~$k\in\N$, and the inequality can be strict. Consider for instance a germ $f\in\End(\C^2,O)$ whose linear term is non-zero but nilpotent; then we have $c(f)=1$ but $c(f^2)\ge 2$. 

Nevertheless, \eqref{eq:2.1} implies the existence of the \emph{asymptotic attraction rate} 
defined by the limit
\[
c_\infty(f)=\lim_{k\to+\infty} c(f^k)^{1/k}\;,
\]
which is a basic (formal and holomorphic, at least) invariant of~$f$. Contrarily to the dimension one case, $c_\infty(f)$ is not necessarily an integer. However, Favre and Jonsson~\cite{FJ1} have
proved the following 

\begin{theorem}[Favre-Jonsson, 2007 \cite{FJ1}]
\label{th:2.FJ1}
Let $f\in\End(\C^2,O)$ be non-invert\-ible and dominant. Then $c_\infty=c_\infty(f)$ is a quadratic integer, i.e., there exist integers $a$,~$b\in\Z$ such that $c^2_\infty+a c_\infty+b=0$. 
Moreover, there exists $\delta\in(0, 1]$ such that $\delta c^k_\infty\le c(f^k)\le c^k_\infty$
for all~$k\ge 1$.
\end{theorem}

This result clearly suggests an open problem:

\begin{itemize}
\item[(OP5)] \emph{Let $f\in\End(\C^n,O)$ be non-invertible and dominant. Is it true that $c_\infty=c_\infty(f)$ is an 
algebraic integer of order at most~$n$, i.e., there are integers $a_0,\ldots, a_{n-1}\in\Z$ such that $c^n_\infty+a_{k-1} c_\infty^{k-1}+\cdots+a_0=0$? Moreover, does there exist 
$\delta\in(0, 1]$ such that $\delta c^k_\infty\le c(f^k)\le c^k_\infty$
for all~$k\ge 1$?}
\end{itemize}

Roughly speaking, the order of a germ may be thought of as a sort of analog of the degree
of a polynomial map: the former somewhat measure the rate of attraction of the origin, while
the latter measure the rate of attraction of infinity. For the sequence of the degree of the
iterates of a polynomial map, Favre and Jonsson \cite{FJ2} have proved the following

\begin{theorem}[Favre-Jonsson, 2011 \cite{FJ2}]
\label{th:2.FJ2}
Let  $F\colon\C^2\to\C^2$ be a polynomial map. Then the sequence $\{\deg F^j\}_{j\in\N}$
satisfies a linear recursion formula with integer coefficients.
\end{theorem}

We can then wonder whether a similar property holds for the sequence of orders of the iterates
of a germ:

\begin{itemize}
\item[(OP6)] \emph{Let $f\in\End(\C^n,O)$ be non-invertible and dominant. Is it true that the sequence 
$\{c(f^k)\}_{k\in\N}$ satisfies, at least for $k$ large enough, a linear recursion formula
with integer coefficients?}
\end{itemize}

\subsection{Non-invertible systems: classification}
\label{sec:2.3}

We have already remarked that in dimension 1 every non-invertible germ of order~$k$ is holomorphically
conjugated to $z^k$, by B\"ottcher's theorem. In several variables, this is not true: 
for instance, as first remarked by Hubbard and Papadopol \cite{HP}, the map
$F(z,w)=(z^2+w^3,w^2)$ cannot be, even topologically, locally conjugated to the homogeneous quadratic
map $H(z,w)=(z^2,w^2)$. Indeed, both maps have as critical locus (which is topologically
defined) the union of the two axes; however the union of the two axes is $H$-invariant but
not $F$-invariant. More precisely, the critical value set of $H$ is the union of the two axes,
whereas the critical value set of $F$ is the union of the $z$-axis with the curve $z^2=w^3$.

\begin{definition}
Given $f\in\End(\mathbb{C}^n,O)$, we shall denote by 
\[
\hbox{\rm Crit}(f)=\{\det(df)=0\}
\] 
the
set of critical points of~$f$, and by
\[
\hbox{\rm PCrit}(f)=\bigcup_{k\ge 0}f^{k}\bigl(\hbox{\rm Crit}(f)\bigr)
\]
the \emph{postcritical set} of~$f$.
\end{definition}

The postcritical set of a homogeneous map is a cone; thus a superattracting germ $f$ can be
topologically (respectively, holomorphically) locally conjugated to a homogeneous map only if
its postcritical set is a topological (respectively, analytic) cone, that is the image of a standard
cone under a local homeomorphism (respectively, biholomorphism). Buff, Epstein and Koch \cite{BEK} have proved that this condition is also sufficient when
the homogeneous map is non-degenerate:

\begin{theorem}[Buff-Epstein-Koch 2011 \cite{BEK}]
\label{th:2.BEK}
Let $f\in\End(\C^n,O)$ be a superattracting germ, and let $H\colon\C^n\to\C^n$ be the first
non-vanishing term, of degree $c=c(f)$, in the homogeneous expansion of~$f$.
Assume that $H$ is non-degenerate, that is $H^{-1}(O)=\{O\}$. Then the following 
assertions are equivalent:
\begin{itemize}
\item[\textup{(i)}] $f$ is holomorphically locally conjugated to~$H$; 
\item[\textup{(ii)}] there is a germ of holomorphic vector field $\xi$ with $\xi(p)=p+o(\|p\|)$ as~$p\to O$
and such that $df\circ\xi=c\, \xi\circ f$; 
\item[\textup{(iii)}] there is a germ of holomorphic vector field $\zeta$ tangent to the postcritical set of~$f$
and such that $\zeta(p)=p+o(\|p\|)$ for~$p\to O$;
\item[\textup{(iv)}] the postcritical set of~$f$ is an analytic cone.
\end{itemize}
\end{theorem}

This result immediately prompts three open problems:

\begin{itemize}
\item[(OP7)] \emph{Let $f\in\End(\C^n,O)$ be a superattracting germ, and let $H$ be the first 
non-zero
term in the homogeneous expansion of~$f$. Assume that $H$ is degenerate, that is there exists
$v\in\C^n\setminus\{O\}$ such that $H(v)=O$. Is it still true that if the postcritical set of~$f$ is
an analytic cone then $f$ is holomorphically locally conjugated to~$H$?}
\end{itemize}

\begin{itemize}
\item[(OP8)]\emph{Let $f\in\End(\C^n,O)$ be a superattracting germ, and let $H$ be the first 
non-zero term in the homogeneous expansion of~$f$.
Assume that the postcritical set of~$f$ is a
topological cone nearby the origin. Under which conditions is then $f$
topologically locally conjugated to~$H$? In particular, is this true when $H$ is non-degenerate?}
\end{itemize}

\begin{itemize}
\item[(OP9)]\emph{Let $f$, $g\in\End(\C^n,O)$ be two superattracting germs, with $c(f)=c(g)$. Is it true that $f$ and $g$ are
topologically (respectively, holomorphically) locally conjugated if and only if their postcritical sets
are topologically (respectively, holomorphically) conjugated, that is, there is a local homeomorphism (respectively, biholomorphism) sending the postcritical set of~$f$ onto the
postcritical set of~$g$?} 
\end{itemize}

Concerning instead the classification problem with respect to the formal conjugacy, \cite{AR}
contains a formal classification of superattracting germs $f\in\End(\C^2,O)$ of
order~2. The methods given there can in principle be used to attack the general formal
classification problem; here we limit ourselves to a more specific question:

\begin{itemize}
\item[(OP10)]\emph{Classify with respect to formal conjugation the superattracting germs
$f\in\End(\C^n,O)$ of order~2 when $n\ge 3$, or of order~$3$ when $n=2$.}
\end{itemize}

\subsection{Non-invertible systems: rigidification}
\label{sec:2.2}

Asking for a holomorphic classification of non-invertible germs is possibly too much;
on the other hand, a birational classification has been obtained by Favre-Jonsson and Ruggiero
at least in dimension~2. Let us introduce a bit of terminology to explain their results.

\begin{definition}
\label{def:2.2.1}
Given $f\in\End(M,p)$, set
\[
\hbox{\rm Crit}^\infty(f)=\bigcup_{k\ge 0}f^{-k}\bigl(\hbox{\rm Crit}(f)\bigr)
=\bigcup_{k\ge 0}\hbox{\rm Crit}(f^k)\;.
\]
We say that $f$ \emph{weakly rigid} if $\hbox{\rm Crit}^\infty(f)$, as a germ at the origin, is composed by $0\le q\le\dim M$ smooth irreducible components, having a simple normal crossing at the origin. If $f$ is weakly rigid and $W_1,\ldots,W_q$ are the irreducible components of 
$\hbox{\rm Crit}^\infty(f)$, we shall say that $f$ is \emph{rigid} if moreover for each $j=1,\ldots,q$ there
exists $I_j\subseteq\{1,\ldots,q\}$ such that $f(W_j)=\bigcap\limits_{i\in I_j} W_i$.
\end{definition}

Invertible germ, because $\textup{Crit}^\infty(f)=\varnothing$, and  
non-dominant germs, because $\textup{Crit}^\infty(f)=M$,  are trivially rigid; so this notion
is interesting only for non-invertible dominant germs.

\begin{definition}
\label{def:2.2.2}
A \emph{modification} of a complex $n$-dimensional manifold $M$ is a surjective 
holomorphic map $\pi\colon\tilde M\to M$ obtained as composition of a finite number 
of blow-ups of submanifolds (or points); in particular, a modification is a birational isomorphism. The modification is \emph{based} at a point~$p_0$
if the first blow-up is made along a submanifold of~$M$ containing~$p_0$, and subsequent
blow-ups are made along submanifolds intersecting the inverse image of~$p_0$. 
\end{definition}

In particular, a modification is a birational isomorphism. In dimension~2 Favre-Jonsson
\cite{FJ1} and Ruggiero \cite{Ru} have proved that every (non-invertible dominant) germ
is birationally conjugated to a rigid germ:
 
\begin{theorem}[Favre-Jonsson, 2007 \cite{FJ1}, Ruggiero, 2009 \cite{Ru}]
\label{th:2.FJR}
For every $f\in\End(\mathbb{C}^2,O)$ there exist a modification $\pi\colon \tilde M\to\C^2$
based at~$O$, a point~$p\in E=\pi^{-1}(O)$ and a 
rigid holomorphic germ $\tilde f\in\End(M,p)$
so that~$\pi\circ\tilde f=f\circ\pi$.
\end{theorem}

Since Favre \cite{Fa} has classified 2-dimensional (non-invertible dominant) rigid germs, this result gives
birational normal forms for non-invertible dominant germs.

Of course, one would like to extend these results at least to dimension~3, and thus we can add 
two more open problems to the list:

\begin{itemize}
\item[(OP11)] \emph{Classify $3$-dimensional non-invertible dominant rigid germs.}
\end{itemize}

\begin{itemize}
\item[(OP12)] \emph{Under which conditions on $f\in\End(\C^3,O)$ there exist a modification $\pi\colon \tilde M\to\C^3$
based at~$O$, a point~$p\in E=\pi^{-1}(O)$ and a 
(possibly weakly) rigid holomorphic germ $\tilde f\in\End(M,p)$
so that~$\pi\circ\tilde f=f\circ\pi$?}
\end{itemize}

It should be mentioned that the latter problem is somewhat related to a non-dynamical result
by Cutkosky \cite{Cu}:

\begin{theorem}[Cutkosky, 2006 \cite{Cu}]
\label{th:2.Cub}
Let $f\colon X\to Y$ be a dominant morphism of algebraic $3$-varieties over~$\C$. Then there
exist: modifications $\phi\colon\tilde X\to X$ and $\psi\colon\tilde Y\to Y$, with $\tilde X$ and
$\tilde Y$ non-singular; simple normal crossing divisors $D_{\tilde Y}$ in~$\tilde Y$ and
$D_{\tilde X}$ in~$\tilde X$; and a morphism $\tilde f\colon \tilde X\to \tilde Y$ such that
the diagram
\[
\xymatrix{\tilde X \ar[r]^{\tilde f} \ar[d]_\phi & \tilde Y \ar[d]^\psi\\
X \ar[r]^f & Y\\}
\]
commutes, $D_{\tilde X}=\tilde f^{-1}(D_{\tilde Y})$ and $\tilde f$ is toroidal with respect to~$D_{\tilde X}$ and~$D_{\tilde Y}$ (i.e., $\tilde f$ is locally given by monomials in suitable \'etale
local parameters on~$\tilde X$).
\end{theorem}

\subsection{Parabolic systems: Fatou flower}
\label{sec:2.4}

Our next topic is the dynamics of \emph{parabolic} systems, where the differential is unipotent;
more precisely we are interested in \emph{tangent to the identity} germs, that is 
holomorphic local dynamical
systems $f\in\End(\mathbb{C}^n,O)$ of the form
\begin{equation}
  \label{eqpsc}
f(z)=z+P_\nu(z)+P_{\nu+1}(z)+\cdots\in\mathbb{C}_0\{z_1,\ldots,z_n\}^n\;,  
\end{equation}
where $P_\nu$ is the first non-zero term in the homogeneous expansion of~$f$,
and $\nu\ge 2$ is the \emph{multiplicity} of~$f$.

Of course, the first thing one would like to do is to find a several variables version
of the Leau-Fatou flower theorem.
To describe the right generalization, we have to introduce a few concepts.

\begin{definition}
  Let $f\in\End(\mathbb{C}^n,O)$ be tangent at the identity and of
  multiplicity~$\nu$. A \emph{characteristic direction} for~$f$ is a non-zero
  vector $v\in\mathbb{C}^n\setminus\{O\}$ such that $P_\nu(v)=\lambda
  v$ for some~$\lambda\in\mathbb{C}$. If $P_\nu(v)=O$ (that is,
  $\lambda=0$) we shall say that $v$ is a \emph{degenerate}
  characteristic direction; otherwise, (that is, if $\lambda\ne 0$) we
  shall say that $v$ is \emph{non-degenerate.} 
\end{definition}

Characteristic directions always exist, and it is not difficult to show (see, e.g., \cite{AT1}) 
that a generic $f$ has exactly  $(\nu^n-1)/(\nu-1)$ characteristic directions, 
counted with respect to a suitable multiplicity. 

The notion of characteristic directions has a dynamical origin. Indeed, it is possible to
prove (see, e.g., \cite{Ha2}) that an orbit 
$\{f^k(z_0)\}$ of a germ $f\in\End(\C^n,O)$ tangent to the identity converges to the origin
\emph{tangentially} to a direction~$[v]\in\mathbb{P}^{n-1}(\mathbb{C})$ --- that is 
$f^k(z_0)\to O$ in~$\mathbb{C}^n$ and $[f^k(z_0)]\to[v]$
in~$\mathbb{P}^{n-1}(\mathbb{C})$, where $[\cdot]\colon\mathbb{C}^n
\setminus\{O\}\to\mathbb{P}^{n-1}(\mathbb{C})$ is the canonical
projection --- then $v$ is a characteristic direction of~$f$.

We can now introduce the several variables analogue of petals: parabolic curves. 

\begin{definition}
   A \emph{parabolic curve} for $f\in\End(\mathbb{C}^n,O)$
tangent to the identity of multiplicity~$\nu\ge 2$ is an injective holomorphic map
$\phe\colon\Delta\to\mathbb{C}^n\setminus\{O\}$ satisfying the following properties:
\begin{itemize}
\item[(a)] \ $\Delta$ is a simply connected domain in~$\mathbb{C}$ with $0\in\partial\Delta$;
\item[(b)] \ $\phe$ is continuous at the origin, and $\phe(0)=O$;
\item[(c)] \ $\phe(\Delta)$ is $f$-invariant, and $(f|_{\phe(\Delta)})^k\to O$
uniformly on compact subsets as $k\to+\infty$.
\end{itemize}
\noindent Furthermore, if $[\phe(\zeta)]\to[v]$ in $\mathbb{P}^{n-1}(\mathbb{C})$ as~$\zeta\to
0$ in~$\Delta$ we shall say that the parabolic curve $\phe$ is \emph{tangent} to
the direction~$[v]\in\mathbb{P}^{n-1}(\mathbb{C})$. 

Finally, a \emph{Fatou flower} tangent to a direction~$[v]$ is a set of $\nu-1$ parabolic curves
tangent to~$[v]$, with domains the connected components of a set of the form
$D_{\delta,\nu}=\{\zeta\in\mathbb{C}\mid
|\zeta^{\nu-1}-\delta|<\delta\}$ for $\delta>0$ small enough.
\end{definition}

Then the first main generalization of the Leau-Fatou flower theorem to several
complex variables is due to \'Ecalle and Hakim (see also \cite{W}):

\begin{theorem}[\'Ecalle, 1985 \cite{E4}; Hakim, 1998 \cite{Ha2}]
\label{th:2.EcalleHakim}  
Let $f\in\End(\mathbb{C}^n,O)$ be a germ
tangent to the identity, and $[v]\in\mathbb{P}^{n-1}(\mathbb{C})$
a non-degenerate characteristic
direction. Then there exists (at least) a Fatou flower for~$f$ tangent to~$[v]$.
\end{theorem}

This result applies to germs tangent to the identity having non-degenerate characteristic directions; however, it is easy to find examples
of germs having only degenerate characteristic directions. In dimension~2 it is possible 
to get Fatou flowers in this case too (see also \cite{BCL}):

\begin{theorem}[Abate, 2001 \cite{A2}]
\label{th:2.Abate}
Every germ $f\in\End(\mathbb{C}^2,O)$
tangent to the identity, with the origin as an isolated fixed point, admits at least
one Fatou flower tangent to some direction.
\end{theorem}

The proof works in dimension~2 only, and this leads to our next open problem:

\begin{itemize}
\item[(OP13)] \emph{Is it true that every germ $f\in\End(\mathbb{C}^n,O)$
tangent to the identity, with the origin as an isolated fixed point, admits at least
one Fatou flower tangent to some direction?}
\end{itemize}

Some partial results are presented in~\cite{AT1} and~\cite{Ro2}, showing however that the analogy with the
local dynamics of holomorphic vector fields that guided the proof of Theorem~\ref{th:2.Abate}
breaks down when $n\ge 3$; so apparently new ideas are needed. 

\subsection{Parabolic systems: classification}
\label{sec:2.5}

Theorems~\ref{th:2.EcalleHakim} and~\ref{th:2.Abate} describe the dynamics only in 1-dimensional subsets of an $n$-dimensional space, and so are very far from determining the
dynamical behavior of a tangent to the identity germ in a full neighbourhood of the origin.

Now, it is possible to attach to each characteristic direction (see \cite{E4} and \cite{Ha3} for
the non-degenerate case, and \cite{ABT1,ABT2,AT3,A7} for the general case)
$n-1$ numbers, called \emph{directors} or \emph{indices} (actually, directors and indices are
not the same numbers, but they are strongly related; see, again, \cite{ABT1} and \cite{AT3}),
that can be useful to describe the behavior in a neighbourhood of the Fatou flowers. For instance,
Hakim~\cite{Ha3} (see also~\cite{ArR}) has proved that if all the directors at a non-degenerate characteristic direction have positive real part then the Fatou flower is attracting, that is there is an open neighbourhood
of the Fatou flower consisting of points whose orbit is converging to the origin tangentially
to the given characteristic direction. On the other hand, there are examples (see, e.g., \cite{AT3} and
references therein) of germs having orbits converging to the origin without being tangent to
any direction, as well as of germs having periodic orbits arbitrarily close to the origin (a
phenomenon that cannot happen in dimension~1). Thus in general the stable set is
larger than the set of points with orbits converging to the origin tangentially to some direction;
however, it the known examples the presence of ``anomalous" points in the stable set 
seems again related to the indices, and in particular to the existence of purely imaginary indices.

Anyway, the natural open problem here is:

\begin{itemize}
\item[(OP14)] \emph{Describe, using characteristic directions, directors, indices and possibly other invariants, the stable set of a germ tangent to the identity in~$\C^n$.}
\end{itemize}

See also \cite{BM} for some results on this problem when all directions are characteristic.

In \cite{AT3} we started a systematic study of the local dynamics of a particularly important class
of dynamical systems tangent to the identity: time-1 maps of homogeneous vector fields. 
Indeed, if we identify, as we can, a homogenous vector field in~$\C^n$ of degree~$\nu\ge 2$ with a $n$-tuple $P_\nu$
of homogeneous polynomials of degree~$\nu$ then its time-1 map is of the form
\[
g(z)=z+P_\nu(z)+O(\|z\|^{\nu+1})\;,
\]
and thus it is tangent to the identity of multiplicity~$\nu$. 

These germs are
particularly important because of the following reformulation of Camacho's Theorem~\ref{th:1.Camacho}:

\begin{corollary}
\label{th:2.Camacho} 
Let $f\in\End(\mathbb{C},0)$ be a holomorphic local dynamical system tangent to the identity
with multiplicity~$\nu$ at the fixed point. Then $f$ is topologically locally
conjugated to the time-$1$ map of the homogeneous vector field
$z^{\nu}\frac{\partial}{\partial z}$.
\end{corollary}

Thus in dimension~one time-1 maps of homogeneous vector fields provide a complete set of
topological normal forms for germs tangent to the identity. This, and the work done in~\cite{AT3},
suggests the following open problem:

\begin{itemize}
\item[(OP15)]\emph{Let $f\in\End(\C^n,O)$ be given by
\begin{equation}
f(z)=z+P_\nu(z)+O(\|z\|^{\nu+1})\;,
\label{eq:2.conj}
\end{equation}
and assume that all characteristic directions of $f$ are non-degenerate. Is it true that $f$ 
is topologically locally conjugated to the time-$1$ map of the homogeneous vector field~$P_\nu$?}
\end{itemize}

The assumption on the characteristic directions is necessary. 
If $v$ is a degenerate characteristic direction for the time-1 map $g$ of a homogeneous
vector field $P_\nu$ (that is, $P_\nu(v)=O$) then the whole complex line $\C v$ consists of
zeroes of the vector field, and thus it is pointwise fixed by~$g$, whereas the fixed point set
of a generic germ of the form (\ref{eq:2.conj}) consists of the origin only.

Concerning the formal classification, in his monumental work \cite{E4} \'Ecalle has given a complete set of formal
invariants for holomorphic local dynamical systems tangent to the identity with
at least one non-degenerate characteristic direction. For instance, he
has proved the following

\begin{theorem}[\'Ecalle, 1985 \cite{E4,E5}]
\label{th:2.Ecalleformal}  
Let $f\in\End(\mathbb{C}^n,O)$ be a holomorphic local dynamical system
tangent to the identity of multiplicity~$\nu\ge 2$. Assume that 
\begin{itemize}
\item[{\rm(a)}] $f$ has exactly $(\nu^n-1)/(\nu-1)$ distinct non-degenerate
characteristic directions and no degenerate characteristic directions;
\item[{\rm(b)}] the directors of any non-degenerate characteristic
direction are irrational and mutually independent over~$\mathbb{Z}$.
\end{itemize}
\noindent Let $[v]\in\mathbb{P}^{n-1}(\mathbb{C})$ be a non-degenerate characteristic
direction, and denote by
$\alpha_1,\ldots,\alpha_{n-1}\in\mathbb{C}$ its directors. Then
there exist a unique $\rho\in\mathbb{C}$ and unique (up to dilations)
formal series~$R_1,\ldots,R_n\in\mathbb{C}[\![z_1,\ldots,z_n]\!]$, where
each $R_j$ contains only monomials of total degree at least~$\nu+1$ and
of partial degree in~$z_j$ at most~$\nu-2$, such that $f$ is formally
conjugated to the time-1 map of the formal vector field
\[
X\!=\!{1\over(\nu-1)(1+\rho z_n^{\nu-1})}\left\{[-z_n^\nu+R_n(z)]\frac{\partial}{\partial
z_n}\!+\!\sum_{j=1}^{n-1}[-\alpha_j z_n^{\nu-1}z_j+R_j(z)]{\partial\over\partial z_j}\right\}.
\]
\end{theorem}

A natural question is how to complete the formal classification:

\begin{itemize}
\item[(OP15)]\emph{Find formal normal forms and invariants for germs tangent to the identity
having only degenerate characteristic directions.}
\end{itemize}

Partial results on this problem can be find in~\cite{BM}, \cite{AT2} and \cite{AR}.

\subsection{Parabolic systems: irregular singularities}
\label{sec:2.6}

In \cite{AT3} we showed that the study of the dynamics of time-1 maps of $n$-dimensional homogeneous vector
fields can be reduced to the study of singular holomorphic foliations in Riemann surfaces 
of~$\PP^{n-1}(\C)$ and of geodesics for meromorphic connections on Riemann surfaces
(notice that a singular holomorphic foliation in $\PP^1(\C)$ is completely determined by its finite
set of singular points, and so when $n=2$ the problem reduces to the study of geodesics for
a meromorphic connection on~$\PP^1(\C)$; this might be the moral reason why the dynamics
of germs tangent to the identity seems to be substantially more
complicated in dimension~3 or more than in dimension~2), and in particular we were able to give a complete
description of the dynamics in a full neighbourhood of the origin for a large class
of 2-dimensional examples. Our study suggested several questions worth of further study;
let us just mention one of them.

We have already remarked that there are two kinds of characteristic directions: non-degenerate and degenerate. Actually, as already noticed in \cite{A2} and exploited in~\cite{ABT1} and~\cite{AT3}, we have to refine this classification because there are different types of degenerate
characteristic directions. To make things simpler, let us assume that $n=2$ and take a homogeneous vector field $Q=(Q_1,Q_2)$ of degree~$\nu\ge 2$.
In \cite{AT3} it is shown that we can reduce the study of the dynamics of $Q$ to the study of
the dynamics of another vector field, the \emph{geodesic field}~$G$, defined on the total space of
a suitable line bundle on $\PP^1(\C)$. In coordinates induced by the usual non-homogeneous coordinates centered in~$[1:0]\in\PP^1(\C)$ the geodesic
field can be written as
\begin{equation}
G= g_1(z)v\frac{\partial}{\partial z}+(\nu-1) g_2(z)v^2\frac{\partial}{\partial v}
\label{eq:2.G}
\end{equation}
where $g_1(z)=Q_2(1,z)-z Q_1(1,z)$ and $g_2(z)=Q_1(1,z)$. In particular, $v_0=(1,\zeta_0)$ is a
characteristic direction of~$Q$ if and only if $g_1(\zeta_0)=0$, and it is a degenerate characteristic 
direction if and only if $g_1(\zeta_0)=g_2(\zeta_0)=0$. 

Let $\mu_j(\zeta_0)\in\N$ be the order of
vanishing of~$g_j$ at~$\zeta_0$. Then we shall say that the point with $(z,v)$ coordinates 
given by $(\zeta_0,0)$ --- which corresponds to the point $[1:\zeta_0]\in\PP^1(\C)$ --- is
\begin{itemize}
\item[--]  an \emph{apparent singularity} if $\mu_1(\zeta_0)\le\mu_2(\zeta_0)$;
\item[--] a \emph{Fuchsian singularity} if $\mu_1(\zeta_0)=\mu_2(\zeta_0)+1$; and
\item[--] an \emph{irregular singularity} if $\mu_1(\zeta_0)>\mu_2(\zeta_0)+1$.
\end{itemize}
Then it is easy to see that non-degenerate characteristic directions are either Fuchsian or 
irregular (whereas degenerate characteristic directions can be apparent, Fuchsian or irregular), and it turns out that Fuchsian singularities with $\mu^1=1$ can be characterized as 
the non-degenerate characteristic directions with non-zero director. 

In \cite{AT3} we found formal normal forms for~$G$ around all kinds of singularities, but holomorphic normal
forms only for apparent and Fuchsian singularities. Using the holomorphic normal forms
we were able to study in detail the local dynamics of $G$ about apparent and Fuchsian singularities
(and thus the dynamics of $Q$ around the corresponding characteristic directions); but it is
still open the study of the local dynamics about irregular singularities. And thus we have our
new open problem:

\begin{itemize}
\item[(OP16)] \emph{Describe the local dynamics of a vector field $G$ of the form \textup{\eqref{eq:2.G}}
about an irregular singularity $(\zeta_0,0)$.}
\end{itemize}

See \cite{Vi} for some results related to this topic.

\begin{remark}
We limited our presentation to germs tangent to the identity; but similar results and problems
can be asked for germs whose differential is represented by a Jordan matrix with 1 as only
eigenvalues. However, in \cite{A1} is shown that such germs are birationally conjugated
(through a modification based at the origin) to germs tangent to the identity, and so many
questions can be reduced to the latter case.

It is also worthwhile to mention that recent works by Bracci, Raissy and Zaitsev \cite{BZ,BRZ}
have shown that techniques developed to study parabolic germs can be also
applied in the study of the dynamics of non-parabolic germs.
\end{remark}

\subsection{Linearizable systems: Brjuno condition}
\label{sec:2.7}

In one variable we saw that every germ whose multiplier was not a root of unity is formally
linearizable. Now, $\lambda\in\C^*$ is a root of unity if and only if $\lambda^q=\lambda$ for
some $q\ge 2$; a similar (but more widespread) phenomenon might prevent formal
linearization in several variables too.

\begin{definition}
\label{def:2.res}
Given $\vec{\lambda}=(\lambda_1,\ldots,\lambda_n)\in\mathbb{C}^n$,
a 
\emph{resonance} for $\vec{\lambda}$ is a relation of the form 
\begin{equation}
\lambda_1^{k_1}\cdots\lambda_n^{k_n}-\lambda_j=0  
  \label{eq:2.cuno}
\end{equation}
for some $1\le j\le n$ and some $k_1,\ldots,k_n\in\mathbb{N}$ with $k_1+\cdots+k_n\ge
2$. If $\lambda_1,\ldots,\lambda_n$ are the eigenvalues of the differential at the origin of
$f\in\End(\mathbb{C}^n,O)$ we shall say that \eqref{eq:2.cuno} is a \emph{resonance of} $f$.
\end{definition}

Resonances are the obstruction to formal linearization. Indeed, a
standard computation shows that the coefficients of a formal
linearization have in the denominators quantities of the form
$\lambda_1^{k_1}\cdots\lambda_n^{k_n}-\lambda_j$; in particular it follows that
a germ $f\in\End(\mathbb{C}^n,O)$ with no
  resonances is always formally conjugated to its
  differential~$df_O$, that is, it is formally linearizable. It should be mentioned that however
a given germ can be formally (and even holomorphically) linearizable even in
presence of resonances; see, e.g., Raissy~\cite{R1,R2}.

The Brjuno problem in several variables consists in deciding when a formal linearization
is actually convergent, keeping in mind that in absence of resonance the formal linearization
is unique, but in presence of resonances 
the formal linearization if it exists is in general not unique; see Raissy~\cite{R3,R4,R5} for a discussion (and much more)
of this problem in a general setting.

To describe the main results known on the Brjuno problem, let us introduce the
following definition:

\begin{definition}
\label{def:2.Bfr}
For $\vec{\lambda}=(\lambda_1,\ldots,\lambda_n)\in\mathbb{C}^n$ and $m\ge 2$ set
\begin{eqnarray}
  \Omega_{\vec{\lambda}}(m)=\min\bigl\{|\lambda_1^{k_1}\cdots
\lambda_n^{k_n}-\lambda_j|&\bigm|& k_1,\ldots,k_n\in\mathbb{N},\,2\le
k_1+\cdots+k_n\le m, \nonumber\\
&& 1\le j\le n,\, \lambda_1^{k_1}\cdots
\lambda_n^{k_n}\ne \lambda_j\bigr\}\;.   
 \label{eq:2.cB}
 \end{eqnarray}
If $\lambda_1,\ldots,\lambda_n$ are the eigenvalues of~$df_O$, we shall
write~$\Omega_f(m)$ for~$\Omega_{\vec{\lambda}}(m)$.
\end{definition}

For some $f\in\End(\C^n,O)$ it might well happen that 
\[
\lim_{m\to+\infty}\Omega_f(m)=0\;,
\]
which is the reason why, even without resonances, the formal
linearization might be diverging, exactly as in the one-dimensional
case. 

Up to now the best result ensuring the convergence of a formal linearization
is the analogue of the Brjuno theorem, that is the implication (iii)$\Longrightarrow$(i)
in Theorem~\ref{th:1.BY}, proved by Brjuno \cite {Brj2,Brj3} in absence of resonances
and generalized by R\"ussmann~\cite{Rus} and Raissy~\cite{R4} to the formally linearizable
case:

\begin{theorem}[Brjuno, 1971 \cite{Brj3}, R\"ussmann, 2002 \cite{Rus}]
\label{Brjunohyp}  
Consider a discrete holomorphic local dynamical system $f\in\End(\mathbb{C}^n,O)$ 
formally linearizable (e.g., without resonances) and with
$df_O$ diagonalizable.  Assume that
\begin{equation}
\sum_{k=0}^{+\infty}{1\over 2^k}\log{1\over\Omega_f(2^{k+1})}<+\infty\;.  
 \label{eq:2.cdue}
\end{equation}
Then $f$ is holomorphically linearizable.
\end{theorem}

\begin{remark}
\label{rem:2.Ynd}
The assumption of diagonalizable differential is necessary. Indeed, Yoccoz \cite{Y2}
has proved that for every $A\in GL(n,\mathbb{C})$ such that its eigenvalues have no resonances and such that its Jordan
normal form contains a non-trivial block associated to an eigenvalue of modulus
one there exists $f\in\End(\mathbb{C}^n,O)$ with $df_O=A$ which is not
holomorphically linearizable.
\end{remark}

Recalling Theorem~\ref{th:1.BY} it is natural to ask whether condition \eqref{eq:2.cdue}
is necessary for the holomorphic linearization of all germs having diagonalizable differential
with given eigenvalues. We can then state the following

\begin{itemize}
\item[(OP17)]\emph{Let $\vec{\lambda}=(\lambda_1,\ldots,\lambda_n)\in\C^n$ be such that
\[
\sum_{k=0}^{+\infty}\frac{1}{2^k}\log\frac{1}{\Omega_{\vec{\lambda}}(2^{k+1})}=
+\infty\;.
\]
Is it possible to find $f\in\End(\C^n,O)$, with diagonalizable differential having eigenvalues
$\lambda_1,\ldots,\lambda_n$, which is formally linearizable but not holomorphically
linearizable?}
\end{itemize}

There are a couple of cases where a positive answer is known. For instance, it is possible
to adapt the classical one-variable construction of Cremer \cite{Cr2} and prove the following

\begin{theorem}
  \label{th:2.BrjunoCremer}
  Let $\vec{\lambda}=(\lambda_1,\ldots,\lambda_n)\in\mathbb{C}$ be without resonances
  and such that
\begin{equation}
\limsup_{m\to+\infty}{1\over
m}\log{1\over\Omega_{\vec{\lambda}}(m)}=+\infty\;.
\label{eq:2.BC}
\end{equation}
Then there exists $f\in\End(\mathbb{C}^n,O)$, with $df_O=\hbox{\rm
diag}(\lambda_1,\ldots,\lambda_n)$, not holomorphically 
linearizable.
\end{theorem}

This does not answer (OP17) because, exactly as in one variable, it is possible to find $\vec{\lambda}\in\C^n$
such that the limsup in \eqref{eq:2.BC} is finite but the series in~\eqref{eq:2.cdue} diverges.

Another easy situation is when one of the components of $\vec{\lambda}$ does not satisfy 
the one-dimensional Brjuno condition. 
Indeed, if $\lambda\in S^1$ does not satisfy \eqref{eqqsei} then any $f\in\End(\mathbb{C}^n,O)$ of
the form
\[
f(z)=\bigl(\lambda z_1+z_1^2, g(z)\bigr)
\]
is not holomorphically linearizable, because if $\phe\in\End(\mathbb{C}^n,O)$ were a holomorphic linearization of~$f$ then $\psi(\zeta)=\phe(\zeta,O)$ would be a holomorphic linearization
of the quadratic polynomial $\lambda\zeta+\zeta^2$, against Theorem~\ref{th:1.BY}.

This again is not enough to solve (OP17) because there are $n$-tuples $\vec{\lambda}\in\C^n$ formed by complex numbers satisfying the one-dimensional Brjuno condition~\eqref{eqqsei} but not satisfying \eqref{eq:2.cdue}: for instance, 
take $\lambda\in S^1$ not a Brjuno number, $0<|\mu|<1$ and put $\lambda_1=\mu$
and $\lambda_2=\lambda\mu^{-1}$. Clearly both $\lambda_1$ and $\lambda_2$ satisfy trivially
\eqref{eqqsei}, whereas 
\[
\lambda_1^{k}\lambda_2^{k+1}-\lambda_2=\mu^{-1}(\lambda^{k+1}-\lambda)\;,
\]
and thus it is easy to see that $(\lambda_1,\lambda_2)$ does not satisfy \eqref{eq:2.cdue}.

We end with a final open problem, which is a recasting in this context of our first open problem:

\begin{itemize}
\item[(OP18)] \emph{Given $\vec{\lambda}=(\lambda_1,\ldots,\lambda_n)\in\C^n$, find families $\{f_{\vec{\lambda},a}\}_{a\in M}\subset\End(\C^n,O)$ of formally 
linearizable germs whose differential is represented by the diagonal matrix of eigenvalues
$(\lambda_1,\ldots,\lambda_n)$ such that 
$f_{\vec{\lambda},a}$ is holomorphically linearizable if and only if $\vec{\lambda}$ satisfies
$\eqref{eq:2.cdue}$.} For instance, \emph{is it true that a germ of the form $f(z)=(\lambda_1 z_1,
\lambda_2 z_2)+Q(z)$, where $Q$ is a pair of quadratic homogeneous polynomials and $(\lambda_1,\lambda_2)$ have no resonances, is holomorphically linearizable if and only if
$(\lambda_1,\lambda_2)$ satisfy~$\eqref{eq:2.cdue}$?}
\end{itemize}



%
%

\begin{acknowledgements}
It is a pleasure to thank Charles Favre, Mattias Jonsson, Jasmin Raissy and Matteo Ruggiero for many useful
discussions and the suggestion of several interesting open problems.
\end{acknowledgements}



\bigskip
{\small\obeylines
\noindent Marco Abate
\noindent Dipartimento di Matematica
\noindent Universit\`a di Pisa
\noindent Largo Pontecorvo 5
\noindent 56127 Pisa
\noindent Italy
\noindent \emph{E-mail:} abate@dm.unipi.it}

\end{document}